\documentclass[12pt]{amsart}

\usepackage{amssymb,amsthm,amsmath,verbatim}
\usepackage{t1enc}
\usepackage{enumerate}
\usepackage{color}
\usepackage{tikz}
\usepackage{breqn}
\usepackage{xmpmulti}
\usepackage{psfrag}

\usepackage{todonotes}

 \usepackage[usenames,dvipsnames]{pstricks}
 \usepackage{epsfig}
\usepackage{pst-grad}
\usepackage{pstricks,pst-node}
\usepackage{float}

\usetikzlibrary{cd}

\newtheorem{theorem}{Theorem}

\newtheorem{claim}[theorem]{Claim}
\newcommand{\dom}{\text{dom}}

\newcommand{\Fn}{{\text Fn}}
\newcommand{\setm}{{\setminus}}

\newcommand{\omg}{{\omega_1}}
\newcommand{\oo}{{\omega}}
\newcommand{\mc}[1]{{\mathcal{#1}}}

\DeclareMathOperator{\pnf}{PONAG_{fin}}
\DeclareMathOperator{\hfc}{HFC}
\DeclareMathOperator{\hfcw}{HFC_w}

\newenvironment{tproof}{
  
  \begin{proof}
}{\end{proof}}

\newenvironment{cproof}{
  
  \begin{proof}
}{\end{proof}}

 \author[D. T. Soukup]{D\'aniel T. Soukup}
\address[D. T. Soukup]{Kurt G\"odel Research Center for Mathematical Logic, Faculty of Mathematics, University of Vienna,  Währinger Strasse 25,
1090 Wien, Austria}
\email[Corresponding author]{daniel.soukup@univie.ac.at}

\urladdr{http://www.logic.univie.ac.at/$\sim  $soukupd73/}

\newtheorem{prob}[theorem]{Problem}
\newtheorem{prop}[theorem]{Proposition}

\author[P. J. Szeptycki]{Paul J. Szeptycki}

\address[P. J. Szeptycki]{Department of Mathematics and Statistics,
Faculty of Science and Engineering, York University, Toronto, Ontario, Canada M3J 1P3}
\email{szeptyck@yorku.ca}
\date{\today}


\begin{document}
\title{A 0-dimensional, Lindel\"of space that is not strongly D}

\begin{abstract}
    A topological space $X$ is strongly $D$ if for any neighbourhood assignment $\{U_x:x\in X\}$, there is a $D\subseteq X$ such that $\{U_x:x\in D\}$ covers $X$ and $D$ is locally finite in the topology generated by $\{U_x:x\in X\}$. We prove that $\diamondsuit$ implies that there is an $\hfcw$ space in $2^\omg$ (hence 0-dimensional, Hausdorff and  hereditarily Lindel\"of) which is not strongly $D$. We also show that any $\hfc$ space $X$ is dually discrete and if additionally countable sets have Menger closure then $X$ is a $D$-space.
\end{abstract}

\maketitle

\section{Introduction} 

A space $X$ is said to be a $D$-space if for every neighbourhood assignment $\{U_x:x\in X\}$ there is a closed discrete set $D\subseteq X$ such that $\{U_x:x\in D\}$ covers the space \cite{vD} (we refer to $D$ as a kernel for the neighbourhood assignment). One of the main open problems  regarding topological covering properties is whether every regular, Lindel\"of space is a $D$-space. The latter question is due to E. van Douwen and we refer the reader to \cite{E,FS,G,HM} for more background.

Recently, L. Aurichi \cite{A} defined a space to be {\em strongly D} if for every neighbourhood assignment $\{U_x:x\in X\}$ there is a set $D\subseteq X$ such that $\{U_x:x\in D\}$ covers $X$ and $D$ is locally finite in the topology generated by $\{U_x:x\in X\}$ i.e., for each $z\in X$, there is a finite $F\subseteq X$ such that 
$$
z\in \bigcap\{U_x:x\in F\} \text{ and } \bigcap\{U_x:x\in F\}\cap D \text{ is finite.}
$$
Note that if the topology generated by $\{U_x:x\in X\}$ is $T_1$ (or if we add the cofinite sets to this basis) then the above condition does imply that $D$ is closed discrete in that topology. Aurichi has shown that every strongly D space is Lindel\"of and that if there is a Lindel\"of, $T_1$ non strongly $D$-space then there is a Lindel\"of, $T_1$ non $D$-space.  On the other hand, in \cite{SSz}, we showed that under the assumption of $\diamondsuit$, there is a $T_2$, hereditarily Lindel\"of non D-space. This provides the closest approximation for a negative solution to van Douwen's question to date. 

\medskip

Now, we present a 0-dimensional, Hausdorff (and hence regular) modification of that example which is still hereditarily Lindel\"of and not strongly $D$.\footnote{Yet another version of this construction was used to show that the union of two $D$-spaces may not be $D$ \cite{SSz2}.} In fact, our space is homeomorphic to an $\hfcw$ subspace of $2^\omg$ which are basic examples of hereditarily Lindel\"of spaces; Section \ref{sec:1} covers all the necessary definitions and the construction itself.

We complement the previous result by showing that any $\hfc$ space (a natural strengthening of being $\hfcw$) is dually discrete i.e., neighbourhood assignments always have discrete kernels (see Theorem \ref{thm:hfcw}). Moreover, we prove that any $\hfc$ space with the property that countable sets have Menger closure is actually a $D$-space (see Theorem \ref{thm:hfc}). The latter results are proved in Section \ref{sec:2}.
Figure \ref{fig:diag} below summarizes the known relations (where HL stands for hereditarily Lindel\"of). 

\begin{figure}[H]

\centering
\begin{tikzcd}
 & \hfc \arrow[d] \arrow[ddrr, blue, bend left, "\textmd{Theorem \ref{thm:hfcw}}"]  & &  \\
  &  \hfcw \arrow[d, blue, "\times"{anchor=center}, "\textmd{Theorem \ref{thm:constr}}\;" left] \arrow[r] & \textmd{HL} \arrow[d, red, "+T_3" near start, "?" description]& \\
 \textmd{Menger}\arrow[r] & \textmd{strongly D} \arrow[r] & \textmd{D-space} \arrow[r] & \textmd{dually discrete}
  \end{tikzcd}
  \bigskip
  
    \caption{The implications between the covering properties}
    \label{fig:diag}
\end{figure}
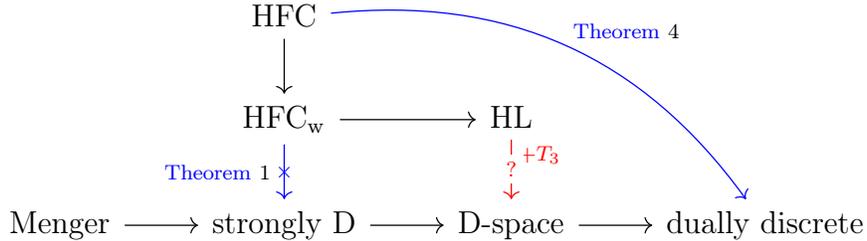

The red arrow marks van Douwen's question and we conclude our paper with a list of further open problems.

\subsection{Notation and terminology} We let $\Fn(I,2)$ denote the set of finite partial function from $I$ to 2. Any such finite function $s$ determines a basic open subset of the Tychonoff product $2^I$ by setting $[s]=\{x\in 2^I:x\supset s\}$.  

For a topological space $(X,\tau)$, a \emph{local $\pi$-network} at $x\in X$ is a family $\mc F$ of arbitrary subsets of $X$ so that for any open neighbourhood $V$ of $x$, there is $F\in \mc F$ with $F\subset V$. A family of open sets $\mc U$ is an \emph{$\oo$-cover} of $X$ if for any finite $F\subset X$, there is $U\in \mc U$ with $F\subset U$.

\subsection{Acknowledgements} D. T. Soukup would like to thank the generous support of the Ontario Trillium Scholarship, FWF Grant I1921 and NKFIH OTKA-113047 during the preparation of this manuscript.

\vskip 13pt
\section{The construction} \label{sec:1}

Our goal is to prove the following.

\begin{theorem}\label{thm:constr} Under $\diamondsuit$, there is a 0-dimensional, Hausdorff and  hereditarily Lindel\"of space which is  not strongly $D$.
\end{theorem}

We will define a topology $\rho$ on $\omg$ by constructing $U_\gamma\subset \omg$ for $\gamma<\omg$. Sets of the form $$
U_s=\left(\bigcap_{\gamma\in s^{-1}(1)} U_\gamma\right) \cap \left(\bigcap_{\gamma\in s^{-1}(0)} \left(\omg\setminus U_\gamma \right)\right)
$$  will form the basis of the topology where $s\in Fn(\omg,2)$. 
Our space $(\omg,\rho)$ naturally embeds into $2^\omg$ by the map  $\beta \mapsto f(\beta):= x_\beta\in 2^\omg$ defined by
\[
x_\beta(\gamma)=
\begin{cases}
1 & \textmd{ if } \beta\in U_\gamma,\\
0 & \textmd{ if } \beta\notin U_\gamma
\end{cases}
\]
where $\gamma<\omg$.
Condition (4) below in our inductive construction ensures that $\rho$ is $T_1$ and in turn, the map $f$ is injective. For any $s\in \Fn(\omg,2)$, the intersection of the basic open set $[s]$ in $2^\omg$ with $f[\omg]$ is exactly $f[U_s]$. So $f$ is a homeomorphism with its image.

To ensure that $(\omg,\rho)$ is hereditarily Lindel\"of, we employ the $\hfcw$ machinery from \cite[Definition 3.2]{J}. Let $\lambda$ be some uncountable cardinal, $F\in[\lambda]^n$ and $b\in 2^n$. We let $F* b$ denote the function from $F$ to 2 which which takes value $b(i)$ on the $i^{th}$ element of $F$. 

Recall that $X\subset 2^\lambda$ is called an $\hfcw$ space if for any $n<\omega$ and uncountable, pairwise disjoint $\mathcal F\subset [\lambda]^n$, there is a countable $\mathcal F_0\subset \mathcal F$ so that for any $b\in 2^n$,

$$|X\setm\bigcup\{[F*b]:F\in \mathcal F_0\}|\leq \omega. $$

Any $\hfcw$ space is hereditarily Lindel\"of \cite[3.3]{J} and these spaces provided some basic combinatorial examples of $L$-spaces (as one of many interesting applications). Let us refer the reader to I. Juh\'asz's \cite{J} for more information on the structure  and properties of such spaces.

We will show in Proposition \ref{prop:hfcw} that $f$ maps $(\omg,\rho)$ to an $\hfcw$ space in $2^\omg$ and so $\rho$ is hereditarily Lindel\"of.


\medskip

The sets $U_\gamma$ will be constructed simultaneously  by an induction of length $\omg$. At each stage $\alpha\geq \gamma$, we will have an approximation $U^\alpha_\gamma$ for $U_\gamma$ and in fact, $U_\gamma\cap (\alpha+1)=U^\alpha_\gamma$. Moreover, we will assume $\gamma+1\subset U_\gamma$ so the subspace topology on $\alpha+1$ will be completely determined by stage $\alpha$.

\medskip

Let us start the construction and assume $\diamondsuit$. Let $\{B_\alpha:\alpha\in \omega_1\}$ be a $\diamondsuit$ sequence capturing subsets of $\Fn(\omega_1,2)$ i.e.,  $B_\alpha\subseteq \Fn(\alpha,2)$ and for any $B\subseteq \Fn(\omega_1,2)$ there are stationary many $\alpha$ such that $B\cap \Fn(\alpha,2)=B_\alpha$. Indeed, we will only be interested in capturing uncountable families $B\subseteq \Fn(\omega_1,2)$ with pairwise disjoint domains in order to assure that the space is $\hfcw$.  
In addition, enumerate all countable subsets of $\omega_1$ as $\{C_\alpha:\alpha\in \omega_1\}$  so that for each $\alpha$, we have that $\sup C_\alpha<\alpha$. These sets are the potentially locally finite kernels that we should avoid.
\medskip

By recursion on $\alpha<\omega_1$, we  define sets $\{U_\gamma^\alpha:\gamma\leq\alpha<\omega_1\}$ so that the following inductive hypotheses are satisfied:
\begin{enumerate}
\item[(1)] For all $\gamma\leq \alpha$, $U_\gamma^\alpha\subseteq \alpha+1$,
\item[(2)] For all $\alpha$ limit and $n<\oo$, $U_{\alpha+n}^{\alpha+n} =(0,\alpha]\cup \{\alpha+n\}$,
\item[(3)] For all $\gamma\leq\alpha<\beta$, $U_\gamma^\beta\cap(\alpha+1)=U_\gamma^\alpha$,
\item[(4)] For all $\eta<\alpha$ there is $\gamma<\alpha$ such that $\eta\in U_\gamma^\alpha$ and $\alpha\not\in U_\gamma^\alpha$ or vice-versa, $\eta\not\in U_\gamma^\alpha$ and $\alpha\in U_\gamma^\alpha$. 
\end{enumerate}
We let $\tau_\alpha$ be the topology on $\alpha+1$ generated by the sets $\{U_\gamma^\alpha:\gamma\leq \alpha\}$ and let $\rho_\alpha$ be the topology on $\alpha+1$ generated by the sets $$\{U_\gamma^\alpha,(\alpha+1)\setm U_\gamma^\alpha:\gamma\leq \alpha\}.$$ Note that by (3), if $\alpha<\beta$ we have that $\tau_\alpha$ is the subspace topology on $\alpha+1$ generated by $\tau_\beta$ and similarly, $\rho_\alpha$ is the subspace topology on $\alpha+1$ generated by $\rho_\beta$. Moreover each $\rho_\alpha$ is zero-dimensional and by (4) also $T_1$.

To present the rest of the inductive hypotheses, we need some more notation and definitions.
For any finite $s\in Fn(\alpha,2)$ let 
$$
U^\alpha_s=\left(\bigcap_{\xi\in s^{-1}(1)} U^\alpha_\xi\right) \cap \left(\bigcap_{\xi\in s^{-1}(0)} \left(\alpha+1\setminus U^\alpha_\xi \right)\right).
$$ 
These sets form a basis for $\rho_\alpha$.
For all $\beta\leq \alpha$, let 
$$
{\mathcal W}_\beta^\alpha = \{U^\alpha_s:s\in B_\beta\}.$$

We will make sure that if $B_\beta$ is large in some sense (see (6) below) then ${\mathcal W}_\beta^\alpha$  covers the interval $(\beta,\alpha]$. This will help us prove hereditarily Lindel\"ofness through the $\hfcw$ machinery.
\medskip

 Now, we have the following additional inductive hypotheses: for all $\beta\leq \alpha$,
\begin{enumerate}
\item[(5)] if $C_\alpha$ is locally finite in  $\tau_\alpha$ then  $\bigcup \{U_\xi^\alpha:\xi\in C_\alpha\}\neq \alpha+1$; 
\item[(6)] if $B_\beta\subseteq Fn(\beta,2)$ consists of functions with pairwise disjoint domains and if there is a countable elementary submodel $M\prec H_{\omega_2}$ that satisfies $M\cap \omega_1=\beta$ and such that 
\begin{itemize}
\item $\{B_\gamma:\gamma<\omega_1\}\in M$,
\item there is an uncountable $B\in M$ such that $B_\beta=B\cap M$, and 
\item there is $(V_\gamma:\gamma<\omega_1)\in M$ such that $V_\gamma\cap \beta= U^\alpha_\gamma \cap \beta$ for all $\gamma<\beta$,
\end{itemize}
 then  
\begin{enumerate}
\item if $\beta \leq \alpha$, then $\{\dom(s):s\in B_\beta\}$ is a local $\pi$-network at $\beta$ in the $\tau_\alpha$ topology, and  
\item if $\beta<\alpha$ then  for each $\tau_\alpha$-neighbourhood $V$ of $\beta$ and each finite subset $F\subseteq V$, 
the following is an $\omega$-cover of $(\beta,\alpha]$:
$$\{U_s^\alpha:s\in B_\beta{\text{ and }}\dom(s)\subseteq V\setminus F\}.$$ 
\end{enumerate}

\end{enumerate}

Let us carry out the construction and verify that the inductive hypotheses can be preserved. First suppose that $\alpha<\omega_1$ and for all $\eta\leq\beta<\alpha$, $U^\beta_\eta$ has been defined satisfying the inductive hypotheses (1)-(6). Let 
$$
{U^{<\alpha}_\eta}=\bigcup\{U^\beta_\eta:\eta<\beta<\alpha\}
$$

and let $\tau_\alpha'$ denote the topology generated by these sets on $\alpha$. Before we continue, note that for each $\beta<\alpha$, $\tau_\beta$ is the subspace topology on $\beta+1$ inherited from $\tau_\alpha'$. Moreover, the key (6)(b) condition about $\oo$-covers  is satisfied by $\tau_\alpha'$ when restricted to the set $(\beta,\alpha)$.

\medskip

Now, for each $\eta<\alpha$, our goal is to extend ${U^{<\alpha}_\eta}$ to $U^\alpha_\eta$ by deciding whether to include or exclude $\alpha$ from it.  Following our previous notation, we let

$$
{U^{<\alpha}_s} = \left(\bigcap_{\xi\in s^{-1}(1)} {U^{<\alpha}_\xi}\right) \cap \left(\bigcap_{\xi\in s^{-1}(0)} \left(\alpha\setminus {U^{<\alpha}_\xi}\right)\right)
$$

for $s\in Fn(\alpha,2)$.

First, let $T$ be the set of $\beta\leq \alpha$ satisfying the hypotheses of (6) and first assume that $\alpha\in T$ and $T\cap \alpha\not=\emptyset$ (hence $\alpha$ is a limit ordinal). Since $\alpha\in T$, let $M$ be the elementary submodel witnessing this and let $B\in M$ be the uncountable family with $B\cap M=B_\alpha$. 

Next, fix an enumeration $\{(\beta_n, G_n) :n\in\omega\}$ of all pairs
$(\beta,G)$ where $\beta\in T\cap\alpha$ and $G$ is a finite subset of the interval $(\beta,\alpha)$ so that each such pair appears infinitely often. 
For each $\beta<\alpha$, we fix a decreasing local neighbourhood base $\{V_n(\beta):n<\omega\}$ in the $\tau_\alpha'$ topology.  If $C_\alpha$ is locally finite in $\alpha$ with respect to $\tau_\alpha'$ then  we may assume that $V_0(\beta)\cap C_\alpha$ is finite (for all $\beta<\alpha$).  

Now to proceed with the construction, we choose a sequence $(s_n)_{n\in \omega}$ by recursion on $n$ such that 
\begin{enumerate}
\item[(i)] $s_n\in B_{\beta_n}$ and $\dom(s_n)\subseteq V_n(\beta_n)$,
\item[(ii)] $\dom(s_n)\cap C_\alpha = \emptyset$ and $ (\sup\bigcup_{k<n}\dom(s_k))+1\notin \dom(s_n)$,
\item[(iii)] $M\models S_{n}:=\{t\in B: \dom(t)\subseteq {U^{<\alpha}_{s_k}}$ for all $k\leq n\}$ is uncountable (and hence the set of $\dom(t)$'s with $t\in B_\alpha$ satisfying the above is infinite), and
\item[(iv)] $G_n\subseteq U^{<\alpha}_{s_n}$.
\end{enumerate} 

Having defined $s_k$ for $k< n$, consider $(\beta_n,G_n)$. By our inductive hypothesis (6)(b), for any $t\in S_{n-1}$ there is some $s\in B_{\beta_n}$ so that 

\begin{equation}\label{eq}
    \dom(s)\subset V_n(\beta_n)\setm (C_\alpha\cup\{(\sup\bigcup_{k<n}\dom(s_k))+1\})
\end{equation} and $\dom(t)\cup G_n\subset U^{<\alpha}_s$. Since $S_{n-1}$ is uncountable, there is a single $s=s_n\in B_{\beta_n}$ as satisfying (\ref{eq}) above, so that $$S_{n}=\{t\in S_{n-1}: \dom (t)\cup G_n\subset U^{<\alpha}_{s_n}\}$$ is uncountable. This argument can be carried out completely in $M$ since all relevant parameters are in $M$ (note that we do not know if $C_\alpha$ is in $M$ but $C_\alpha \cap V_n(\beta_n)$ is finite so an element of $M$). This defines $s_n$ satisfying (i)-(iv) and concludes this induction.


Now we are ready to define $U^\alpha_\eta$ for all $\eta\leq \alpha$. Let

\[
U^\alpha_\eta:=
\begin{cases}
U^{<\alpha}_\eta\cup {\alpha}, & \textmd{ if } \eta\in \bigcup_{n\in \oo}s_n^{-1}(1),\\
U^{<\alpha}_\eta, & \textmd{ if } \eta\in \alpha\setm \bigcup_{n\in \oo}s_n^{-1}(1),\\
\alpha+1, & \textmd{ if } \eta=\alpha.
\end{cases}
\]




We must check that the inductive hypotheses (1)-(6) are satisfied. Items (1)-(3) follow directly from the construction. Regarding (4), since $\eta\in U^\alpha_\beta$ for all $\eta<\beta$, if the set of $\beta$ such that $\alpha\not\in U^\alpha_\beta$ is cofinal in $\alpha$ then (4) would be satisfied. Hence item (4) follows as long as the sets $s_n^{-1}(1)$ do not cover a tail of $\alpha$. This is made sure by condition (ii) in the construction of $s_n$.


Item (5) has been taken care of in the construction of $s_n$ by condition (ii): indeed if $C_\alpha$ is $\tau_\alpha$ locally finite then it is so in $\tau_\alpha'$ as well. Then $\dom(s_n)\cap C_\alpha=\emptyset$ and so $\alpha\not\in U_\xi^\alpha$ for all $\xi\in C_\alpha$.


To verify (6) first consider (6)(a). For $\beta<\alpha$, this is ensured by the inductive hypothesis and the fact that $\tau_\alpha$ and $\tau_\beta$ neighbourhoods coincide on $\beta+1$. So we consider $\alpha$ and fix a $\tau_\alpha$ neighbourhood $V$ of $\alpha$. There is a finite set $F\subseteq \bigcup\{s_k^{-1}(1):k<\omega\}$ such that $U_F=\bigcap_{\xi\in F} U^\alpha_\xi\subseteq V$. Therefore, there is an $n$ such that $\bigcap_{k\leq n}U^\alpha_{s_k}\subseteq U_F$. By (iv) in the construction of the $s_n$, we have many $t\in B_\alpha$ such that $\dom(t)\subseteq \bigcap_{k\leq n}U^\alpha_{s_k}\subseteq U_F\subseteq V$. Hence (6)(a) is satisfied.

Finally, to check 6(b), suppose $\beta<\alpha$ and $\beta$ satisfies the hypotheses of the statement. In turn, $\beta\in T$. Fix a $\tau_\alpha$ neighbourhood $V$ of $\beta$ and let $G\subseteq (\beta,\alpha]$ be finite. Since the pair $(\beta,G\setminus\{\alpha\})$ was enumerated infinitely often in the construction, there is an $n$ such that $(\beta,G\setminus\{\alpha\})=(\beta_n,G_n)$ and such that $V_n(\beta)\subseteq V$. Therefore, the chosen $s_n\in B_\beta$ satisfies that $\dom(s_n)\subseteq V$ and $$G_n\subseteq U^{<\alpha}_s\subseteq U^\alpha_s.$$ Looking at the definition of $U^\alpha_\eta$ for $\eta\in \dom(s_n)$, we see that $\alpha\in U^\alpha_\eta$ iff $s(\eta)=1$. This translates to $\alpha\in U^\alpha_{s_n}$ and so $G\subset U^\alpha_{s_n}$ as desired. Hence 6(b) is verified and the main inductive step is finished.


\medskip

Previously, we dealt with the case when $\alpha\in T$ and $T\cap \alpha\not=\emptyset$ where $T$ was the set of $\beta\leq \alpha$ which satisfied the assumptions of condition (6). Now, if $\alpha\not\in T$ but $\alpha$ is a limit and $T\cap \alpha\not=\emptyset$ then the construction is  more simple. We define the sequence of $s_n$ similarly but without  arranging that the domains of elements of $B_\alpha$ form a $\pi$-network at $\alpha$ (i.e., condition (iii) is not needed). 

If $T\cap \alpha=\emptyset$ then the construction is even simpler. Nothing needs to be done except declare $U_\alpha^\alpha$ according to (2) and declare that $ U^\alpha_\eta=U^{<\alpha}_\eta$ for all $\eta<\alpha$. The successor case is handled trivially as well using condition (2). 
\medskip

This ends the construction. The final topology $\rho$ on $\omg$ is generated by the sets $$U_\gamma=\bigcup_{\gamma\leq \alpha<\omg}U^\alpha_\gamma$$ and their complements for $\gamma<\omg$. To reiterate the beginning of the proof, for any $s\in \Fn(\omg,2)$, we have a basic clopen set 
 $$
U_s=\left(\bigcap_{\gamma\in s^{-1}(1)} U_\gamma\right) \cap \left(\bigcap_{\gamma\in s^{-1}(0)} \left(\omg\setminus U_\gamma \right)\right).
$$ 

The topology is clearly  0-dimensional and also $T_1$ by condition (4). We are left to show that $\rho$ is hereditarily Lindel\"of and not strongly $D$.

\begin{prop}\label{prop:hfcw}For any uncountable family $B\subseteq \Fn(\omg,2)$ with pairwise disjoint domains, there is a countable $B'\subset B$ so that $$\omg\setm \bigcup\{U_F:F\in B'\}$$ is countable.

\end{prop}

Now the $\hfcw$ property is easily verified for the homeomorphic subspace $f[\omg]$ of $ 2^\omg$ and so the topology $\rho$ is hereditarily Lindel\"of.

\begin{proof} Given $B$, we can find an elementary submodel $M\prec H(\Theta)$ so that $M$ contains all the relevant parameters and for 
$\beta=M\cap \omg$, $B\cap \Fn(\beta,2)=B_\beta$. Indeed, there are club many models with all the parameters and $B$ is guessed stationary often by the $\diamondsuit$ sequence. In turn, at every stage $\alpha\geq \beta$, $\beta$ satisfied the assumptions of condition (6) and so we made sure that $\bigcup\{U^\alpha_s:s\in B_
\beta\}$ covers $(\beta,\alpha]$ (see (6)(b)). Hence, $$\omg\setm \bigcup\{U_s:s\in B_
\beta\}\subset \omg\setm (\beta+1),$$
as desired.






\end{proof}

\begin{prop}
The topology $\rho$ is not strongly $D$.
\end{prop}
\begin{cproof}
Indeed, this is witnessed by the neighbourhood assignment $\alpha\mapsto U_\alpha$. Suppose that $C$ is locally finite in the topology generated by $\{U_\alpha:\alpha<\omg\}$. Since $\rho$ is Lindel\"of and $C$ is locally finite in $\rho$ too, $C$ must be countable. Hence, there is an $\alpha<\omg$ so that $C=C_\alpha$. In turn, at step $\alpha$ of the main induction, we made sure that $\alpha\notin U^\alpha_\xi=U_\xi\cap (\alpha+1)$ for $\xi\in C_\alpha$ (since $C=C_\alpha$ was $\tau_\alpha$ locally finite at that point). So $X\neq \bigcup_{\xi\in C_\alpha}U_\alpha$, as desired.
\end{cproof}

This concludes the proof of Theorem \ref{thm:constr}.
We remark that any countable subspace of our topology $\rho$ is second countable. In turn, our space is dually second countable i.e., any neighbourhood assignment has a second countable kernel.

\section{$\hfc$ and $D$-spaces}\label{sec:2}

Our goal in this section is to analyse a strengthening of the $\hfcw$ property: a subspace $X\subseteq 2^\lambda$ is called \emph{$\hfc$} if for every $n\in \omega, b\in 2^n$ and any infinite family  $\mathcal F\subset [\lambda]^{n}$ of  pairwise disjoint sets,  $$X\setminus\bigcup \{[F*b]:F\in \mathcal F\}$$ is countable \cite[Definition 3.2]{J}. 
Any $\hfc$ space is $\hfcw$ and so  hereditarily Lindel\"of as well.


While $\hfc$ spaces are not necessarily left separated, in many cases we get left separated spaces in the classical constructions of hereditarily Lindel\"of spaces. It is well-known that every left-separated space is a $D$-space \cite{G}. In general, we do not know whether all $\hfc$ spaces are $D$-spaces. However, we have the following result.

\begin{theorem}\label{thm:hfcw}
Any $\hfc$ space is dually discrete.
\end{theorem}
Recall that $X$ is \emph{dually discrete} if for any neighbourhood assignment $\{U_x:x\in X\}$ there is a discrete $D\subseteq X$ so that $\{U_x:x\in D\}$ covers $X$. It is also unknown whether all (hereditarily) Lindel\"of space are dually discrete \cite{O}.

\begin{tproof}
Given an $\hfc$ space $X\subset 2^\lambda$, we can assume that  neighbourhood assignments are of the form $N:X\to \Fn(\lambda,2)$. We start by taking a countable elementary submodel $M\prec H(\Theta)$ (with $\Theta$ appropriately large) so that $X,N,\lambda\in M$. Our first goal is to find a discrete $D\subset M\cap X$ so that $X\setm N[D]$ is countable. 

List all pairs $(\varepsilon,b)\in Fn(\lambda\cap M,2)\times 2^{<\omega}$ as $\{(\varepsilon_{2\ell},b_{2\ell}):\ell\in \omega\}$ such that each $(\varepsilon,b)$ pair appears infinitely often. We construct a sequence of points $(x_n)_{n<\omega}$ in $M\cap X$  so that
\begin{enumerate}
	\item\label{it:disc} $x_n\in X\cap M\setm \bigcup\{[N(x_k)]:k<n\}$, 
	\item $M\cap X\subset \bigcup\{[N(x_n)]:n\in \oo\}$ and
\item\label{it:meet}	if $n$ is even and there is $x\in X\setm \{x_k:k<n\}$ so that 
	\begin{enumerate}[(i)]
	    \item $\varepsilon_n\subseteq N(x)$,
	    \item $F=\dom(N(x)\setm \varepsilon_n)$ is disjoint from $\dom(N(x_k))$ for $k<n$, and 
	    \item  $N(x)\setm \varepsilon_n=F* b_n$
	\end{enumerate} then we choose $x_n$ to be such.

%
\end{enumerate} 
	
	Note that condition (\ref{it:disc}) ensures that $D=\{x_n:n<\oo\}$ is discrete. The construction is simple: at odd stages we work towards covering $M\cap X$ and at even stages, we see if condition (\ref{it:meet}) can be satisfies: if so, we pick such an $x_n$, otherwise an arbitrary one.
	
	
	We use $N[D]$ to denote $\bigcup\{[N(x)]:x\in D\}$.
	
\begin{claim}\label{clm:remainder0} $|X\setm N[D]|\leq \omega$.
\end{claim}


\begin{cproof}
Suppose otherwise. Then there is $Z\in [X\setm N[D]]^{\omega_1}$ and $(\varepsilon,b)\in Fn(\lambda\cap M)\times 2^{<\omega}$ such that for all  $x\in Z$,

\begin{enumerate}
    \item $N(x)\cap M=\varepsilon$, and 
    \item $N(x)\setm \varepsilon=\dom(N(x)\setm \varepsilon)* b$.
\end{enumerate}

Let $\Gamma=\{n\in \omega: (\varepsilon,b)=(\varepsilon_n,b_n)\}$ and recall the inductive construction of the sequence $\{x_n:n\in\omega\}$. 
In particular, the set $Z$ witnesses that we were able to choose $x_n$ according to condition (\ref{it:meet})  when $n\in \Gamma$.

In turn, there is an infinite set $\tilde D\subseteq D$ so that $\varepsilon\subseteq N(x)$ for each $x\in\tilde D$, $\{\dom(N(x)\setm \varepsilon):x\in \tilde D\}$ is pairwise disjoint and $N(x)\setm \varepsilon=\dom(N(x)\setm \varepsilon)* b$ for all $x\in \tilde D$.

Now, since $X$ is an $\hfc$ space, $$X\setm \bigcup\{[\dom(N(x)\setm \varepsilon)* b]:x\in \tilde D\}$$ is countable. In particular, there is $z\in Z$ such that $z\in \bigcup\{[\dom(N(x)\setminus \varepsilon)* b]:x\in \tilde D\}$. Pick $x\in \tilde D$ such that $z\in [\dom(N(x)\setminus \varepsilon)* b]$. Recall that $z\in [N(z)]\subset [\varepsilon]$ and hence $z\in [\varepsilon]\cap [\dom(N(x)\setminus \varepsilon)* b]=[N(x)]$. This contradicts $z\in Z\subseteq X\setm N[D]$.

\end{cproof}

	Now, list $X\setm N[D]$ as $\{z_n:n\in\oo\}$ (if $N[D]$ already covers $X$ then the proof is done). We define $y_n\in M\cap X$  so that
	\begin{enumerate}
	\setcounter{enumi}{2}
	    \item $z_n\in [N(y_n)]$ and
	    \item $y_n\notin \bigcup_{\ell\leq k_n} [N(x_\ell)]$
	\end{enumerate} where $k_n$ is the maximum of $n$ and $\min\{k<\oo:y_{n-1}\in őN(x_k)]\}$.
	 Why is this possible? At step $n$, we consider the family of open sets $$\{[N(y)]:y\in X\setm \bigcup_{\ell\leq k_{n}} [N(x_\ell)]\}.$$ Note that the latter is in $M$ and  since $X$ is hereditarily Lindel\"of, there is a countable subfamily in $M$ with the same cover. In turn, we can pick $y=y_n\in X\cap M\setm \bigcup_{\ell\leq k_n} [N(x_\ell)]$ which covers $z_n$.
	 
	 \begin{claim}
	 $\{x_n,y_n:n<\oo\}$ is discrete.
	 \end{claim}
	 \begin{cproof}
	 Simply note that $\bigcup_{\ell\leq k_{n+1}} [N(x_\ell)]$ is a neighbourhood of both $x_n$ and $y_n$ which contains only finitely many other $x_k,y_k$.
	 \end{cproof}
	
	This finishes the proof of the theorem since $$X=\bigcup\{[N(x_n)],[N(y_n)]:n<\oo\}.$$
	
\end{tproof}

Our final theorem shows that in a class of  'locally small' topologies, any  $\hfc$ space must be a $D$-space.

\begin{theorem}\label{thm:hfc} Suppose that $X\subset 2^\lambda$ is  \emph{$\hfc$} and the closure of every countable subset of $X$ is Menger.\footnote{The latter property is a natural variant of the well-studied {\em $\omega$-boundedness} assumption i.e., countable sets having compact closure. See \cite{JMW} for a fairly recent overview. We thank L. Zdomskyy for recommending this version of our result; our original theorem assumed that countable sets in $X$ have $\sigma$-compact closure.} Then $X$ is a $D$-space.
\end{theorem}

Recall that a space $Y$ is {\em Menger} if for any countable sequence of open covers $(\mc U_n)_{n\in \oo}$, there are finite $\mc V_n\subseteq \mc U_n$ such that $\bigcup_{n\in \oo}\mc V_n $ covers $Y$. Let us refer the interested reader to \cite{Sch} for background in selection principles and topology.

Any $\sigma$-compact space or Lindel\"of space of size $<\mathfrak d$ is Menger.\footnote{Here $\mathfrak d$ denotes the dominating number \cite{B}.} Moreover, Aurichi showed that all Menger spaces are $D$-spaces \cite[Corollary 2.7 ]{A0} and we need some of the topological games that he used in the proof. The {\em partial open neighbourhood assignment game} (or PONAG, for short) is played as follows. Our two players are $\textbf{N}$ and $\textbf{C}$ and $\textbf{N}$ starts by playing a partial neighbourhood assignment $\{V_x:x\in Y_0\}$ for some $Y_0\subset X$  covering $X$. Next, $\textbf{C}$ replies by a $D_0\subset Y_0$  closed discrete in $X$. In general, $\textbf{N}$ plays $\{V_x:x\in Y_n\}$ for some  $Y_n\subset X\setm \bigcup\{V_x:x\in D_k,k<n\}$ so that  $\{V_x:x\in Y_n\}$  covers $ X\setm \bigcup\{V_x:x\in D_k,k<n\}$ and $\textbf{C}$ replies by $D_n\subset Y_n$ closed discrete in $X$. 
\medskip

\begin{center}
\begin{tabular}{c|c c c c c c}
     Player $\textbf{N}$ & $\{V_x:x\in Y_0\}$ &  & $\{V_x:x\in Y_1\}$ & \dots \\
     \hline 
     Player $\textbf{C}$ &                   & $D_0\subset Y_0$ &  &  $D_1\subset Y_1$ & \dots
\end{tabular}
\end{center}
\medskip

Player $\textbf{C}$ wins if $\bigcup\{V_x:x\in D_n,n<\oo\}$ covers $X$.

We will use the fact that if $Y$ is Menger then $\textbf{N}$ has no winning strategy in PONAG \cite[Proposition 2.6]{A0}; from this, $Y$ being a $D$-space easily follows. In fact, Aurichi's proves in \cite[Proposition 2.6]{A0} that $\textbf{N}$ has no winning strategy in the following modification $\pnf$ of the original PONAG game: $\textbf{N}$ plays as before but $\textbf{C}$ is only allowed to reply by finite sets (instead of arbitrary closed discrete ones). The winning condition is the same as before. This minor modification is quite important in our following proof (and also shows that Menger spaces are strongly $D$).

\begin{tproof}[Proof of Theorem \ref{thm:hfc}] 
First, any neighbourhood assignment (after some shrinking) can be coded by a map  $N:X\to \Fn(\lambda,2)$ so that $x\in [N(x)]$. We use $N[E]$ to denote the set $\bigcup \{[N(x)]:x\in E\}$ in short.

Our plan is to find a closed discrete set $D$ such that $X\setm N[D]$ is countable. Since countable spaces are $D$-spaces, we can find a closed discrete $\tilde D\subset X\setm N[D]$ so that $N[D\cup \tilde D]=X$, as desired (note that the union $ D\cup \tilde D$ is still closed discrete).

Now, let $M\prec H(\Theta)$ be a countable elementary submodel (with $\Theta$ appropriately large) so that $X,N,\lambda\in M$. Let $Y=\overline{X\cap M}$ and note that $Y$ is Menger.


We construct a sequence $(\hat D_n)_{n<\omega}$ of finite sets in $M\cap X$  so that 
\begin{enumerate}
	\item $Y\subset N[\hat D]$ for $\hat D=\bigcup\{\hat D_n:n\in \omega\}$, and
	\item $\hat D_n\cap N[\hat D_{n-1}]=\emptyset$.
 	\end{enumerate} 
We will find these sets using the fact that player $\textbf{N}$ has no winning strategy in the $\pnf$ game on $Y$. So, we define a strategy $\sigma$ for $\textbf{N}$ as follows. We set $V_x=[N(x)]$ and let $\textbf{N}$ start by playing a cover $(V_x)_{x\in Y_0}$ of $X$ for some countable $Y_0\in M$. This is possible since $X$ is Lindel\"of (and we can choose the witness for that in $M$ by elementarity) and so $Y_0\subset M$ as well. In general, given a partial game play $Y_0,\hat D_0,Y_1,\hat D_1,\dots,Y_{n-1},\hat D_{n-1}$ so that $Y_k\in M$, $\sigma$ will simply play $(V_x)_{x\in Y_n}$ where $Y_n\subset X\setm \bigcup_{k<n}N[\hat D_k]$ is a countable set in $M$ so that $(V_x)_{x\in Y_n}$ covers $X\setm \bigcup_{k<n}N[\hat D_k]$. This is again possible by the Lindel\"of property and elementarity.\footnote{This is the point where we need the modified $\pnf$ game so that the sets $\hat D_n$ are in $M$ too.} Now, we know this strategy cannot be winning so there is a run  $Y_0,\hat  D_0,Y_1, \hat D_1,\dots$ following the strategy $\sigma$ so that player $\textbf{C}$ wins i.e., $Y\subset N[\hat D]$ for $\hat D=\bigcup\{\hat D_n:n\in \omega\}$.

Now, we shall enlarge $\hat D$ carefully to cover the rest of $X$ (modulo a countable set). List all pairs $(\varepsilon,b)\in Fn(\lambda\cap M,2)\times 2^{<\omega}$ as $\{(\varepsilon_n,b_n):n\in \omega\}$ such that each such pair appears infinitely often. By induction on $n<\oo$, we define $D_n\supseteq \hat D_n$ by adding at most one point $x_n$ from $M\cap X\setm (\bigcup_{k<n} N[D_k]\cup N[\hat D_{n}])$. How is $x_n$ selected? We look at $(\varepsilon_n,b_n)$ and see if there is $x\in X\setm (\bigcup_{k<n} N[D_k]\cup N[\hat D_{n}])$ so that 
	\begin{enumerate}[(i)]
	    \item $\varepsilon_n\subseteq N(x)$,
	    \item $F=\dom(N(x)\setm \varepsilon_n)$ is disjoint from $\dom(N(y))$ for all $y\in \bigcup_{k<n} D_k\cup \hat D_{n}$, and 
	    \item  $N(x)\setm \varepsilon_n=F* b_n$.
	\end{enumerate} If yes, then we pick such an $x_n$ in $M$.

\medskip



We shall prove that $D=\bigcup\{D_n:n\in \omega\}$ is as required.

\begin{claim} $D$ is locally finite and in turn, closed and discrete.
\end{claim}
\begin{cproof}
Since $D\subseteq M$, it suffices to show that $D$ is locally finite in $Y=\overline{X\cap M}$. Let $y\in Y$ and pick $n\in \omega$ such that $y\in N[\hat D_n]$. Then $U=N[\hat D_n]$ is an open neighbourhood of $y$ and $U\cap D\subseteq \bigcup_{k\leq n} D_k$ is finite.
\end{cproof}

\begin{claim}\label{clm:remainder} $|X\setm N[D]|\leq \omega$.
\end{claim}

In fact, the open sets corresponding to the extra points $x_n$ we chose cover $X$ modulo a countable set.\footnote{So in the end, the Menger machinery is used to make this sequence closed discrete and the $\hfc$ property to cover $X$.} The proof of this claim is exactly as the proof of Claim \ref{clm:remainder0} which we omit repeating.

This concludes the proof of the theorem.
\end{tproof}

We do not know how much the assumption of $X$ being $\hfc$ can be weakened or if the local smallness assumption can be dropped.

\section{Open problems}

The main problem of van Douwen remains open.

\begin{prob}[van Douwen]
Is there a regular, Lindel\"of but non $D$-space? 
\end{prob}

Surprisingly, the following question seems to be open as well (see \cite{O}).

\begin{prob}
Is there a $T_1$, hereditarily Lindel\"of but non dually discrete space? 
\end{prob}

We should emphasise that even consistent examples would be very welcome. 
\medskip

Regarding Theorem \ref{thm:hfc} and Theorem \ref{thm:hfcw}, we ask the following.

\begin{prob}
Is every $\hfc$ space a $D$-space?
\end{prob}

\begin{prob}
Is every $\hfcw$ space dually discrete?
\end{prob}

\begin{prob}
Suppose that $X$ is an $\hfc$ space and every countable subset of $X$ has Menger closure. Is $X$ strongly $D$?
\end{prob}

\begin{prob}
Suppose that every countable subset of $X$ has Menger closure. Is $X$ a $D$-space whenever it is
\begin{enumerate}[(a)]
    \item (hereditarily) Lindel\"of, or 
    \item $\hfcw$?
\end{enumerate}
\end{prob}

It would be equally interesting to see an answer for the preceding questions if instead of the Menger property we assume that countable sets have $\sigma$-compact closure.
Finally, let us refer the reader to \cite{G} for many more interesting open problems around $D$-spaces.


\begin{thebibliography}{10}
\bibitem{O} Alas, Ofelia T., Lucia R. Junqueira, and Richard G. Wilson. {\em Dually discrete spaces} Topology and its Applications 155.13 (2008): 1420-1425.
\bibitem{A0} Aurichi, Leandro F. {\em D-spaces, topological games, and selection principles} Topology Proc. Vol. 36. 2010.

\bibitem{A} L.F.\ Aurichi, {\em D-spaces, separation axioms and covering properties}

\bibitem{B} Blass, Andreas. {\em Combinatorial cardinal characteristics of the continuum} Handbook of set theory. Springer, Dordrecht, 2010. 395-489.


Houston Journal of Mathematics, {\bf 37} No. 3 (2011) 1035-1042.
\bibitem{vD} E.K.\ van Douwen and W.F.\ Pfeffer, {\em Some properties of the Sorgenfrey line and related spaces} Pacific Journal of Mathematics, {\bf 81} No.2 (1979) 371-377.
\bibitem{Do} A.\ Dow {\em An introduction to applications of elementary submodels to topology} Topology Proc. 13 (1988), no. 1, 17--72.
\bibitem{E} T.\ Eisworth {\em On D-spaces} in Open Problems in Topology, Elsevier, (2007) 129-134. 
\bibitem{FS} W.G.\ Fleissner and A.M.\ Stanley {\em D-spaces} Topology and its Applications 114 (2001) 261--271.
\bibitem{G} G.\ Gruenhage {\em A survey of D-spaces}, Contemporary Mathematics, to appear. 
\bibitem{HM} M.\ Hru\v s\'ak and J.T. Moore {\em Introduction: Twenty problems in set-theoretic topology} in Open Problems in Topology, Elsevier, (2007) 111-113. 

\bibitem{J} Juh\'asz, Istv\'an. {\em {HFD} and {HFC} type spaces, with applications} Topology and its Applications 126.1-2 (2002): 217-262.

\bibitem{JMW} Juh\'asz, Istv\'an, Jan van Mill, and William Weiss. {\em Variations on $\omega$-boundedness } Israel Journal of Mathematics 194.2 (2013): 745-766.

\bibitem{K} K.\ Kunen {\em Set Theory, An Introduction to Independence Proofs} Studies in Logic and the Foundations of Mathematics v102, North-Holland, 1983.
\bibitem{vMTW} J.\ van Mill, V.V.\ Tkachuk and R.G. Wilson {\em Classes defined by stars and neighbourhood assignments} Topology Appl. 154 (2007), 2127--2134.

\bibitem{Sch} Scheepers, Marion. {\em Selection principles and covering properties in topology} Note di Matematica 22.2 (2003): 3-41.


\bibitem{SSz} Soukup, D\'aniel T., and Paul J. Szeptycki. {\em A counterexample in the theory of D-spaces} Topology and its Applications 159.10-11 (2012): 2669--2678.

\bibitem{SSz2}  Soukup, D\'aniel T., and Paul J. Szeptycki. {\em The union of two D-spaces need not be D} Fundamenta Mathematicae 220.2 (2013): 129-137.


\end{thebibliography}
\end{document}